\documentclass[11pt]{article}
\usepackage{amsmath,amssymb,amsthm}
\usepackage{mathrsfs}
\usepackage[margin=1in]{geometry}

\theoremstyle{plain}
\newtheorem{theorem}{Theorem}[section]
\newtheorem{lemma}{Lemma}[section]
\newtheorem{proposition}{Proposition}[section]
\newtheorem{corollary}{Corollary}[section]

\theoremstyle{definition}
\newtheorem{definition}{Definition}[section]
\newtheorem{remark}{Remark}[section]
\newtheorem{example}{Example}[section]

\newcommand{\dom}{\operatorname{dom}}
\newcommand{\lsc}{lower semicontinuous}

\begin{document}

\title{Summable Orbits and the Minimal Caristi Potential}
\author{Robledo Mak's Miranda Sette\\
	Faculdade de Ci\^encias Exatas e Tecnologia\\
	Universidade Federal da Grande Dourados -- UFGD\\
	Rodovia Itahum, km 11, Dourados--MS, Brazil\\
	\texttt{robledosette@ufgd.edu.br}}
\date{}
\maketitle

\begin{abstract}
Let $(M,d)$ be a complete metric space and $f:M\to M$. Associated with $f$ is the
\emph{orbit potential}
\[
\varphi_f(x)=\sum_{n\ge0} d\bigl(f^n(x),f^{n+1}(x)\bigr)\in[0,+\infty],
\]
whose finiteness at a single point expresses the summability of the corresponding
forward orbit. We show that, whenever $\varphi_f$ is \lsc, the map $f$ has a fixed point
if and only if some orbit is summable. The lower semicontinuity of each individual gap
$x\mapsto d(f^n(x),f^{n+1}(x))$ is a convenient sufficient condition for this hypothesis,
and we exhibit an example in which it fails while $\varphi_f$ remains \lsc, so that the
criterion applies strictly beyond that condition. Under the same hypothesis we observe
that the existence of a summable orbit is equivalent to $f$ being a Caristi map, and that
$\varphi_f$ is then the \emph{minimal} Caristi potential, in the sense that
$\varphi_f\le\varphi-\inf\varphi$ for every admissible potential $\varphi$. Finally we
delimit the reach of the criterion among generalized contractions: it recovers the
Bianchini--Grandolfi contractions (those governed by a summable comparison function), and
in particular the Banach contraction principle, but it does not subsume the Boyd--Wong or
Matkowski classes, whose orbits need not be summable; the comparison function
$\psi(t)=t/(1+t)$ marks this boundary explicitly.
\end{abstract}


\section{Introduction}

Metric fixed point theory rests on a few structural theorems that connect the existence of
solutions of $f(x)=x$ with the geometry of the underlying space. The Banach contraction
principle~\cite{Banach1922} provides a constructive existence and uniqueness criterion in
complete metric spaces, together with convergence of the Picard iterates. A second, more
flexible principle emerged from variational analysis: Caristi's fixed point
theorem~\cite{Caristi1976}, which replaces contractivity by the existence of a scalar
potential controlling the displacement $d(x,f(x))$. Caristi's theorem is equivalent to
Ekeland's variational principle~\cite{Ekeland1974,Ekeland1979} and to Takahashi's
minimization principle~\cite{Takahashi1991}; these equivalences place fixed point existence
inside the framework of variational methods (see~\cite{GoebelKirk,KhamsiKirk,Rus2001} for
systematic accounts).

In practice the potential in Caristi's theorem may be difficult to produce, while a purely
dynamical quantity is often directly computable: the total displacement along a forward
orbit,
\begin{equation}\label{eq:orbit-sum}
\varphi_f(x)=\sum_{n\ge0} d\bigl(f^n(x),f^{n+1}(x)\bigr).
\end{equation}
When \eqref{eq:orbit-sum} is finite for some $x$, the orbit $(f^n(x))_n$ is Cauchy, and in a
complete space it converges. Whether the limit is a fixed point is precisely the question
addressed here. The idea that a summable orbit forces a fixed point is classical and
underlies both the proof of the Banach principle and the theory of Picard
operators~\cite{Rus2001,Berinde2007}; our aim is not to claim it as new but to isolate the
minimal hypothesis under which it holds, to relate it cleanly to Caristi's theorem, and to
determine exactly which generalized contractions it covers.

The paper is organized as follows. Section~\ref{sec:prelim} records two elementary facts
about the orbit potential: an additive identity along orbits and a semicontinuity lemma.
Section~\ref{sec:main} states the orbit-summability criterion under the hypothesis that
$\varphi_f$ is \lsc\ and shows, by an example, that this is strictly weaker than requiring
each orbit gap to be \lsc. Section~\ref{sec:caristi} identifies the equivalence with
Caristi's theorem and shows that $\varphi_f$ is the minimal Caristi potential.
Section~\ref{sec:applications} recovers the Bianchini--Grandolfi contractions, hence the
Banach principle, and marks the boundary of the method with an explicit example. Concluding
remarks are collected in Section~\ref{sec:remarks}.

\section{Preliminaries}\label{sec:prelim}

Throughout, $(M,d)$ is a metric space and $f:M\to M$ is an arbitrary map. We write
$f^0=\mathrm{id}_M$ and $f^{n+1}=f\circ f^n$. For $n\ge0$ define the \emph{$n$-th orbit gap}
\[
d_f^n(x)=d\bigl(f^n(x),f^{n+1}(x)\bigr),\qquad x\in M,
\]
and the \emph{orbit potential} $\varphi_f:M\to[0,+\infty]$ by \eqref{eq:orbit-sum}.
Recall that a function $\varphi:M\to(-\infty,+\infty]$ is \emph{proper} if
$\varphi(x)<+\infty$ for some $x$, and \emph{\lsc} if
$\varphi(x)\le\liminf_{y\to x}\varphi(y)$ for every $x$.

\begin{lemma}[Orbit identity]\label{lem:cocycle}
For every $x\in M$,
\begin{equation}\label{eq:cocycle}
\varphi_f(x)=d(x,f(x))+\varphi_f(f(x)),
\end{equation}
as an identity in $[0,+\infty]$. Consequently, for every $n\ge0$,
$\varphi_f(f^n(x))=\sum_{j\ge n}d_f^j(x)$.
\end{lemma}

\begin{proof}
Re-indexing the series of nonnegative terms,
\[
\varphi_f(f(x))=\sum_{n\ge0}d\bigl(f^n(f(x)),f^{n+1}(f(x))\bigr)
=\sum_{n\ge0}d\bigl(f^{n+1}(x),f^{n+2}(x)\bigr)
=\sum_{m\ge1}d_f^m(x).
\]
Adding the term $m=0$, namely $d_f^0(x)=d(x,f(x))$, yields \eqref{eq:cocycle}. Iterating
\eqref{eq:cocycle} gives $\varphi_f(f^n(x))=\sum_{j\ge n}d_f^j(x)$.
\end{proof}

\begin{lemma}[Semicontinuity of the orbit potential]\label{lem:lsc}
If each orbit gap $d_f^n$ is \lsc, then $\varphi_f$ is \lsc.
\end{lemma}

\begin{proof}
The partial sums $\varphi_f^N=\sum_{n=0}^N d_f^n$ are finite sums of \lsc\ functions, hence
\lsc. Since the terms are nonnegative, $(\varphi_f^N)_N$ is nondecreasing and
$\varphi_f=\sup_N\varphi_f^N$. The pointwise supremum of any family of \lsc\ functions is
\lsc\ (its epigraph is the intersection of closed epigraphs), so $\varphi_f$ is \lsc.
\end{proof}

\begin{remark}
Lemma~\ref{lem:lsc} is only a sufficient condition; the criterion of the next section requires
merely that $\varphi_f$ itself be \lsc. Example~\ref{ex:separator} below shows that this is
genuinely weaker.
\end{remark}

\section{The orbit-summability criterion}\label{sec:main}

\begin{definition}
A point $x\in M$ has a \emph{summable orbit} if $\varphi_f(x)<+\infty$, that is, if
\[
\sum_{n\ge0}d\bigl(f^n(x),f^{n+1}(x)\bigr)<+\infty.
\]
\end{definition}

\begin{theorem}[Orbit-summability criterion]\label{thm:main}
Let $(M,d)$ be complete and $f:M\to M$ a map such that $\varphi_f$ is \lsc. Then $f$ has a
fixed point if and only if some point of $M$ has a summable orbit.
\end{theorem}

\begin{proof}
($\Leftarrow$) Suppose $\varphi_f(x_0)<+\infty$ and set $x_n=f^n(x_0)$. For $m>n$,
\[
d(x_n,x_m)\le\sum_{j=n}^{m-1}d(x_j,x_{j+1})\le\sum_{j\ge n}d_f^j(x_0).
\]
As $\varphi_f(x_0)=\sum_{j\ge0}d_f^j(x_0)<+\infty$, the tail $\sum_{j\ge n}d_f^j(x_0)\to0$,
so $(x_n)$ is Cauchy; by completeness $x_n\to p$ for some $p\in M$. By
Lemma~\ref{lem:cocycle}, $\varphi_f(x_n)=\sum_{j\ge n}d_f^j(x_0)\to0$. Since $\varphi_f$ is
\lsc,
\[
0\le\varphi_f(p)\le\liminf_{n\to\infty}\varphi_f(x_n)=0,
\]
hence $\varphi_f(p)=0$. As $\varphi_f(p)$ is a sum of nonnegative terms, its first term
vanishes: $d(p,f(p))=d_f^0(p)=0$, i.e.\ $f(p)=p$.

($\Rightarrow$) If $f(p)=p$, then $d_f^n(p)=d(p,p)=0$ for all $n$, so
$\varphi_f(p)=0<+\infty$; the orbit of $p$ is summable.
\end{proof}

\begin{remark}
The forward implication is trivial (a fixed point is its own summable orbit); the content of
the theorem is the reverse one. Only the lower semicontinuity of $\varphi_f$ at limits of orbits is used,
together with completeness and the identity \eqref{eq:cocycle}.
\end{remark}

The next example shows that the hypothesis ``$\varphi_f$ is \lsc'' is strictly more general
than ``each $d_f^n$ is \lsc'' (Lemma~\ref{lem:lsc}).

\begin{example}\label{ex:separator}
Let
\[
M=\{0\}\cup\Bigl\{\tfrac1n:n\ge1\Bigr\}\cup\Bigl\{\tfrac1n+4^{-n}:n\ge1\Bigr\}\cup\{2\}
\subset\mathbb{R},
\]
with the usual metric. The only accumulation point of $M$ is $0\in M$, so $M$ is closed in
$\mathbb{R}$, hence complete; every point except $0$ is isolated. Define $f:M\to M$ by
\[
f(2)=2,\qquad f(0)=2,\qquad f\bigl(\tfrac1n\bigr)=\tfrac1n+4^{-n},\qquad
f\bigl(\tfrac1n+4^{-n}\bigr)=2.
\]
The gap $d_f^0(x)=d(x,f(x))$ satisfies $d_f^0(0)=2$ but
$d_f^0\bigl(\tfrac1n\bigr)=4^{-n}\to0$ while $\tfrac1n\to0$; hence
$\liminf_{x\to0}d_f^0(x)=0<2=d_f^0(0)$, so $d_f^0$ is \emph{not} \lsc\ at $0$. On the other
hand, the orbit of every point terminates at the fixed point $2$, and a direct computation
gives $\varphi_f(0)=2$ and, for $n\ge1$,
\[
\varphi_f\bigl(\tfrac1n\bigr)=4^{-n}+\bigl(2-\tfrac1n-4^{-n}\bigr)=2-\tfrac1n,\qquad
\varphi_f\bigl(\tfrac1n+4^{-n}\bigr)=2-\tfrac1n-4^{-n}.
\]
Both expressions tend to $2=\varphi_f(0)$ as $n\to\infty$, and all other points are isolated,
so $\varphi_f$ is \lsc\ (indeed continuous). Thus Theorem~\ref{thm:main} applies and yields a
fixed point (namely $2$), even though the sufficient condition of Lemma~\ref{lem:lsc} fails.
\end{example}

\section{Relation to Caristi's theorem}\label{sec:caristi}

We use Caristi's condition in \emph{descent form}, which avoids ambiguous differences at
points where the potential is infinite.

\begin{definition}
A map $f:M\to M$ is a \emph{Caristi map} if there exists a proper, bounded-below, \lsc\
function $\varphi:M\to(-\infty,+\infty]$ such that
\begin{equation}\label{eq:caristi}
\varphi(f(x))\le\varphi(x)-d(x,f(x))\qquad(x\in M).
\end{equation}
Such a $\varphi$ is called a \emph{Caristi potential} for $f$. On $\dom\varphi=\{\varphi<+\infty\}$,
\eqref{eq:caristi} coincides with the usual form $d(x,f(x))\le\varphi(x)-\varphi(f(x))$.
\end{definition}

\begin{theorem}[Caristi~\cite{Caristi1976}]\label{thm:caristi}
Every Caristi map on a complete metric space has a fixed point.
\end{theorem}

Our two elementary lemmas turn the orbit potential into the canonical Caristi potential.

\begin{proposition}[Summable orbit $\Rightarrow$ Caristi]\label{prop:os-to-caristi}
Suppose $\varphi_f$ is \lsc\ and some point has a summable orbit. Then $\varphi_f$ is a
Caristi potential for $f$.
\end{proposition}

\begin{proof}
By hypothesis $\varphi_f$ is proper and \lsc, and it is bounded below by $0$. For the descent
inequality \eqref{eq:caristi} with $\varphi=\varphi_f$, use the identity
\eqref{eq:cocycle}: if $\varphi_f(x)<+\infty$ then
$\varphi_f(f(x))=\varphi_f(x)-d(x,f(x))$, an equality; if $\varphi_f(x)=+\infty$ then the
right-hand side of \eqref{eq:caristi} is $+\infty\ge\varphi_f(f(x))$. In either case
\eqref{eq:caristi} holds.
\end{proof}

\begin{proposition}[Caristi $\Rightarrow$ summable orbit; minimality]\label{prop:caristi-to-os}
If $\varphi$ is any Caristi potential for $f$, then for every $x\in M$
\begin{equation}\label{eq:minimal}
\varphi_f(x)\le\varphi(x)-\inf_M\varphi.
\end{equation}
In particular every $x\in\dom\varphi$ has a summable orbit, and $\varphi_f$ is the minimal
Caristi potential in the sense of \eqref{eq:minimal}.
\end{proposition}

\begin{proof}
Write $m=\inf_M\varphi>-\infty$. If $\varphi(x)=+\infty$, \eqref{eq:minimal} is trivial. Let
$\varphi(x)<+\infty$. Applying \eqref{eq:caristi} to $f^j(x)$ and summing over
$j=0,\dots,N$ telescopes to
\[
\sum_{j=0}^N d_f^j(x)\le\varphi(x)-\varphi(f^{N+1}(x))\le\varphi(x)-m.
\]
Letting $N\to\infty$ gives $\varphi_f(x)\le\varphi(x)-m<+\infty$.
\end{proof}

Combining Propositions~\ref{prop:os-to-caristi} and~\ref{prop:caristi-to-os} with
Theorem~\ref{thm:main} yields the following equivalence.

\begin{theorem}\label{thm:equiv}
Let $(M,d)$ be complete and $f:M\to M$ with $\varphi_f$ \lsc. Then the following are
equivalent:
\begin{enumerate}
\item[\textup{(i)}] some point of $M$ has a summable orbit;
\item[\textup{(ii)}] $f$ is a Caristi map.
\end{enumerate}
When they hold, $\varphi_f$ is a Caristi potential for $f$, it is minimal in the sense of
\eqref{eq:minimal}, and $f$ has a fixed point.
\end{theorem}

\begin{proof}
(i)$\Rightarrow$(ii) is Proposition~\ref{prop:os-to-caristi} (with potential $\varphi_f$).
(ii)$\Rightarrow$(i) is Proposition~\ref{prop:caristi-to-os}. Minimality is
\eqref{eq:minimal}, and the existence of a fixed point follows from
Theorem~\ref{thm:main} (or from Theorem~\ref{thm:caristi}).
\end{proof}

\begin{remark}
Theorem~\ref{thm:equiv} makes precise the informal statement that, under semicontinuity of the
orbit potential, ``being a Caristi map'' and ``having a summable orbit'' are the same property,
and that the smallest potential witnessing it is the orbit sum itself. The identification of
$\varphi_f$ as the natural Lyapunov function of a Picard iteration is well known
(cf.~\cite{Rus2001,Berinde2007}); \eqref{eq:minimal} records its extremal position among all
Caristi potentials.
\end{remark}

\section{Application: summable comparison contractions}\label{sec:applications}

We now determine which generalized contractions fall under Theorem~\ref{thm:main}. The
relevant class is governed by comparison functions whose iterates are summable.

\begin{definition}[\cite{Bianchini1968,Berinde2007}]
A function $\psi:[0,\infty)\to[0,\infty)$ is a \emph{summable comparison function} (a
$(c)$-comparison function, or Bianchini--Grandolfi gauge) if $\psi$ is nondecreasing,
$\psi(0)=0$, $\psi$ is continuous at $0$, and $\sum_{n\ge0}\psi^n(t)<+\infty$ for every
$t\ge0$.
\end{definition}

For such $\psi$ one has $\psi(t)<t$ for all $t>0$: otherwise $\psi(t)\ge t$ together with
monotonicity gives $\psi^n(t)\ge t$ for all $n$, contradicting $\psi^n(t)\to0$.

\begin{corollary}[Bianchini--Grandolfi contractions]\label{cor:bg}
Let $(M,d)$ be complete and let $f:M\to M$ satisfy
\begin{equation}\label{eq:bg}
d(f(x),f(y))\le\psi\bigl(d(x,y)\bigr)\qquad(x,y\in M)
\end{equation}
for a summable comparison function $\psi$. Then $f$ has a unique fixed point, and every
Picard orbit converges to it.
\end{corollary}

\begin{proof}
Since $\psi(0)=0$ and $\psi$ is continuous at $0$, \eqref{eq:bg} makes $f$ continuous; hence
each $d_f^n$ is continuous, and $\varphi_f$ is \lsc\ by Lemma~\ref{lem:lsc}. Fix $x\in M$ and
put $c_n=d_f^n(x)$. From \eqref{eq:bg} applied to the pair $(f^n(x),f^{n+1}(x))$,
\[
c_{n+1}=d\bigl(f^{n+1}(x),f^{n+2}(x)\bigr)\le\psi\bigl(d(f^n(x),f^{n+1}(x))\bigr)=\psi(c_n).
\]
Since $\psi$ is nondecreasing, induction gives $c_n\le\psi^n(c_0)$, whence
\[
\varphi_f(x)=\sum_{n\ge0}c_n\le\sum_{n\ge0}\psi^n(c_0)<+\infty.
\]
Thus $x$ has a summable orbit, and Theorem~\ref{thm:main} provides a fixed point $p$, which is
the limit of $(f^n(x))_n$. For uniqueness, if $f(p)=p$ and $f(q)=q$ with $t=d(p,q)>0$, then
$t=d(f(p),f(q))\le\psi(t)<t$, a contradiction; hence $p=q$.
\end{proof}

\begin{corollary}[Banach contraction principle~\cite{Banach1922}]\label{cor:banach}
Every map $f$ on a complete metric space with $d(f(x),f(y))\le c\,d(x,y)$ for a fixed
$c\in[0,1)$ has a unique fixed point.
\end{corollary}

\begin{proof}
Take $\psi(t)=c\,t$, so that $\sum_{n\ge0}\psi^n(t)=t/(1-c)<+\infty$; apply
Corollary~\ref{cor:bg}.
\end{proof}

The following example marks the exact boundary of the method: not every generalized
contraction has summable orbits, so Theorem~\ref{thm:main} does not subsume the Boyd--Wong
or Matkowski classes.

\begin{example}\label{ex:boundary}
Let $\psi(t)=\dfrac{t}{1+t}$. Then $\psi$ is continuous, increasing, and $\psi(t)<t$ for
$t>0$, so it is a Boyd--Wong gauge~\cite{BoydWong1969} and a Matkowski comparison
function~\cite{Matkowski1975} (indeed $\psi^n(t)=t/(1+nt)\to0$). Hence any $f$ satisfying
$d(f(x),f(y))\le\psi(d(x,y))$ on a complete space has a unique fixed point by the theorems of
Boyd--Wong or Matkowski. However,
\[
\psi^n(t)=\frac{t}{1+nt}\sim\frac1n\qquad(n\to\infty),
\]
so $\sum_{n\ge0}\psi^n(t)=+\infty$ for every $t>0$: the corresponding orbit bound is
\emph{not} summable, and $\psi$ is not a summable comparison function. Thus these classes lie
outside the scope of Theorem~\ref{thm:main}, which reaches exactly the contractions governed
by a summable comparison rate.
\end{example}

\section{Concluding remarks}\label{sec:remarks}

The orbit-summability criterion (Theorem~\ref{thm:main}) isolates the summability of a single
forward orbit as the operative hypothesis for the existence of a fixed point, the only
regularity required being lower semicontinuity of the orbit potential. Theorem~\ref{thm:equiv}
shows that, under this regularity, orbit summability and Caristi's condition describe the same
maps, with the orbit sum as the minimal potential. Corollary~\ref{cor:bg} places the
Bianchini--Grandolfi contractions -- hence the Banach principle -- within this framework,
while Example~\ref{ex:boundary} shows that comparison functions with non-summable iterates,
such as $\psi(t)=t/(1+t)$, escape it. The criterion is therefore a faithful but bounded tool:
it certifies fixed points precisely when the geometric displacement of an orbit is finite.

Two directions are left open. First, since Caristi's theorem characterizes metric
completeness~\cite{Weston1977,Kirk1976}, one may ask whether the equivalence of
Theorem~\ref{thm:equiv} does likewise. Second, the identity \eqref{eq:cocycle} is specific to
the iteration of a single map; a nonautonomous version, applicable to sequences generated by
varying maps (as in relaxed or averaged iterations), would require replacing $\varphi_f$ by a
cocycle adapted to the schedule and is not treated here.

\end{document}